\documentclass[11pt]{amsart}
\usepackage{amscd}
\usepackage{amssymb}
\usepackage{hyperref}

\newtheorem{thm}{Theorem}[section]
\newtheorem{lem}[thm]{Lemma}

\newtheorem{cor}[thm]{Corollary}

\newtheorem{prop}[thm]{Proposition}

\theoremstyle{definition}

\renewcommand{\thecase}{}

\newtheorem{conj}[thm]{Conjecture}

\newtheorem{defn}[thm]{Definition}

\renewcommand{\thestep}{}

\theoremstyle{remark}

\makeatletter
\def\alphenumi{
  \def\theenumi{\alph{enumi}}
  \def\p@enumi{\theenumi}
  \def\labelenumi{(\@alph\c@enumi)}}
\makeatother

\makeatletter
\def\thecase{\@arabic\c@case}
\makeatother

\numberwithin{equation}{section}

\makeatletter
\def\thestep{\@arabic\c@step}
\makeatother

\newcommand\barmu{{\bar\mu}}

\newcommand\AAA{\mathbb{A}}

\newcommand\CC{\mathbb{C}}
\newcommand\DD{\mathbb{D}}
\newcommand\EE{\mathbb{E}}
\newcommand\FF{\mathbb{F}}

\newcommand\RR{\mathbb{R}}

\newcommand\ZZ{\mathbb{Z}}

\newcommand\bS{{\mathbf{S}}}

\newcommand\thalf{{\textstyle{\frac{1}{2}}}}

\newcommand\tquarter{{\textstyle{\frac{1}{4}}}}

\newcommand\fg{{\mathfrak{g}}}

\newcommand\fs{{\mathfrak{s}}}

\newcommand\ft{{\mathfrak{t}}}

\newcommand\la{\lambda}
\newcommand\La{\Lambda}
\newcommand\ka{\kappa}

\newcommand\Om{\Omega}
\newcommand\si{\sigma}
\newcommand\Si{\Sigma}

\newcommand\su{{\mathfrak{s}\mathfrak{u}}}

\newcommand\SO{\operatorname{SO}}

\newcommand\U{\operatorname{U}}

\newcommand\ch{\operatorname{ch}}

\newcommand\Coker{\operatorname{Coker}}

\newcommand\Hom{\operatorname{Hom}}

\newcommand\Ker{\operatorname{Ker}}

\newcommand\PD{\operatorname{PD}}

\newcommand\Stab{\operatorname{Stab}}

\newcommand\Sym{\operatorname{Sym}}

\newcommand\Tor{\operatorname{Tor}}

\newcommand\id{{\mathrm{id}}}

\newcommand\spinc{\text{$\text{spin}^c$ }}
\newcommand\spinu{\text{$\text{spin}^u$ }}

\newcommand\Spinu{\text{$\text{Spin}^u$}}

\newcommand\sA{{\mathcal{A}}}
\newcommand\sB{{\mathcal{B}}}

\newcommand\sE{{\mathcal{E}}}
\newcommand\sF{{\mathcal{F}}}
\newcommand\sG{{\mathcal{G}}}

\newcommand\sN{{\mathcal{N}}}

\newcommand\sV{{\mathcal{V}}}

\newcommand\tM{{\tilde M}}

\input{xy}
\xyoption{all}

\begin{document}
\title[Degeneracy loci of families of Dirac operators]
{Degeneracy loci of families of Dirac operators}
\author[Thomas G. Leness]{Thomas G. Leness}
\address{Department of Mathematics\\
Florida International University\\
Miami, FL 33199}
\email{lenesst@fiu.edu}
\urladdr{http://www.fiu.edu/$\sim$lenesst}
\dedicatory{}
\thanks{The author was supported in part by a Florida International University
Faculty Research Grant.}
\date{This version: \today. math.DG/yymmnnn.}
\begin{abstract}
Generalizing some results from R. Leung's thesis, \cite{RLeungThesis},
we compute, in rational cohomology, the Poincare dual of the degeneracy locus
of the family of Dirac operators parameterized by
the moduli space of projectively anti-self-dual $\SO(3)$ connections.
This is the first step in a program to derive a relation between the
Donaldson and spin invariants.
\end{abstract}

\date{April 17th, 2008}
\maketitle

\section{Introduction}
The definition of the spin invariants of a smooth four-manifold,
due to V. Pidstrigach and A. Tyurin in \cite{PTDirac}, is sufficiently
similar to that of the Donaldson invariant to suggest that
there must a formula relating the two.
In this note, we perform a computation which will be useful in deriving
such a relation. This relation, together with
relations between the Donaldson and Seiberg-Witten invariants
and between the spin and Seiberg-Witten invariants given by the $\SO(3)$
monopole cobordism, \cite{FL2a,FL2b,FL5}, should be an important tool for
deriving more explicit forms of these relations and possibly  constraints
on these invariants.

For $E^w_\ka\to X$  a complex, rank-two vector bundle over
a smooth, closed, oriented, four-manifold,  the Donaldson invariants,
\cite{DonPoly,KMStructure}, are
defined by integrating $\mu$-classes over the moduli space $M^w_\ka$ of
projectively anti-self-dual connections on $E^w_\ka$.
If $\fs=(\rho,W)$ is a \spinc structure on $X$, then each unitary connection
$A$ on $E^w_\ka$ defines a Dirac operator on the \spinu structure
$\ft=(\rho\otimes\id_E,W\otimes E)$ of index $n_a(\ft)$.  If $n_a(\ft)\le 0$,
then by \cite{FeehanGenericMetric,TelemanGenericMetric}
the subspace  $J^{\La,w}_\ka\subset M^w_\ka$,
defined in \ref{eq:DefineJ}, of connections whose  Dirac operator has
one-dimensional kernel is,
for generic perturbations, a smooth submanifold of
codimension $2(1-n_a(\ft))$. In this note,  we generalize
some results of \cite{RLeungThesis} by
computing the Poincare dual of  $J^{\La,w}_\ka$.

\begin{thm}
\label{thm:Main}
Let $E^w_\ka\to X$ be a complex, rank-two vector bundle over a
smooth, closed, oriented, four-manifold with $b_1(X)=0$.
Let $M^w_\ka$ be the moduli space of projectively anti-self-dual connections on $E^w_\ka$
and let $K^w_\ka\subset M^w_\ka$ be any compact, codimension-zero submanifold.
Let $\ft$ be a \spinu structure on $X$ with $p_1(\ft)=\ka$, $c_1(\ft)=\La$,
and $n_a(\ft)\le 0$.
Let $J^{\La,w}_\ka\subset M^w_\ka$ be the degeneracy locus of the
\spinu structure $\ft$.  Then, as an element of rational cohomology,
the Poincare dual of $J^{\La,w}_\ka$ in $K^w_\ka$ is
\begin{equation}
(-1)^{1-n_a(\ft)}\sum_{i+2j+2k=1-n_a(\ft)} f_{i,2j,2k} \mu(\ft)^i\smile \Om^j\smile \mu(x)^k,
\end{equation}
where $\mu(\ft)$ and $\Om$  are the $\mu$-classes defined in
\eqref{eq:DefineMuTandOm}
and the coefficients $f_{i,j,k}$ are determined by the recursion relation
\eqref{eq:Recursion1}
and are given as
$$
F(x,y,z)=\sum_{i,j,k} f_{i,2j,2k} x^iy^{2j}z^{2k}
=\exp(\thalf xJ_1(z)+\tquarter y^2J_2(z)+J_3(z)),
$$
where the functions $J_i(z)$ are given by
\begin{align*}
J_1(z){}&= z^{-1}\tan^{-1}(z)
\\
J_2(z){}&=z^{-3}\left( z-\tan^{-1}(z)\right)
\\
J_3(z){}&=-\thalf n_a(\ft)\ln(1+z^2)+\ka(z^{-1}\tan^{-1}(z)-1).
\end{align*}
\end{thm}

The spin invariants of \cite{PTDirac} are defined by integrating $\mu$-classes
over the Uhlenbeck closure of $J^{\La,w}_\ka$:
$$
P^{\La,w}_X(z)=\langle \barmu(z),[\bar J^{\La,w}_\ka]\rangle.
$$
Theorem \ref{thm:Main} implies that the Poincare dual of $J^{\La,w}_\ka$ can
be written as an expression in $\mu$-classes, which we denote as $\mu(T^{\La,w}_\ka)$.  One might
thus expect the Donaldson and spin invariants to be related by
\begin{equation}
\label{eq:Incorrect}
P^{\La,w}_X(z)=\langle \barmu(z),[\bar J^{\La,w}_\ka]\rangle
=\langle \barmu(z)\smile\barmu(T^{\La,w}_\ka),[\bar M^w_\ka]\rangle
=
D^w_X(z T^{\La,w}_\ka).
\end{equation}
However, R. Leung's thesis, \cite{RLeungThesis},
shows that \eqref{eq:Incorrect} only holds when $n_a(\ft)=0$.
If $\sV(T^{\La,w}_\ka)$ is the geometric representative of $\mu(T^{\La,w}_\ka)$
used to define the Donaldson invariant, then Theorem \ref{thm:Main} shows
that the intersections of $\sV(T^{\La,w}_\ka)$ and $J^{\La,w}_\ka$ with
any compact subset of $M^w_\ka$ are cobordant but the same is not
true for their closures in the Uhlenbeck compactification.
Thus, while Theorem \ref{thm:Main} does not give a relation between
the Donaldson and spin invariants, it does allow one to localize
the error in \eqref{eq:Incorrect} near the lower strata of the Uhlenbeck
compactification. In \S \ref{sec:Compare}, we discuss the role this localization
plays in producing a formula for the error in \eqref{eq:Incorrect}.

We do not expect such a formula for the error to be explicit but rather one
of the type appearing in the Kotschick-Morgan conjecture, \cite{KotschickMorgan},
or in the $\SO(3)$-monopole cobordism formula, \cite{FL5}.  The
$\SO(3)$-monopole cobordism formula gives relations between both the
Donaldson and Seiberg-Witten invariants and between the spin and Seiberg-Witten
invariants but these relations are given by unknown polynomials with
universal coefficients depending on topological invariants of $X$.
These coefficients have been determined in some cases but not all.
Having a direct relation between the Donaldson and spin invariants would
give us additional leverage for determining these coefficients.
Indeed, it is possible that this additional leverage will reveal new, topological
constraints of the type appearing in \cite{FKLM} on these invariants of the smooth structure.

\section{Preliminaries}
\label{sec:prelim}

\subsection{ASD moduli space}
Let $E^w_\ka\to X$ be a complex, rank-two, Hermitian vector bundle
over a smooth, closed, oriented four-manifold with $b^1(X)=0$,
$w=c_1(E^w_\ka)$ and $\ka=\langle c_2(E^w_\ka)-\tquarter c_1(E^w_\ka)^2,[X]\rangle$.
Fix a connection $A^w$ on the line bundle $\det(E^w_\ka)$.
For $k\ge 2$, define $\sA^w_\ka$ to be the $L^2_k$
completion of the space of unitary connections $A$ on
$E^w_\ka$ with $A^{\det}=A^w$.  Let $\sA^{w,*}_\ka\subset \sA^w_\ka$
be the irreducible connections.
Define $\sG^w_\ka$ to be the
$L^2_{k+1}$ completion of the group of
special-unitary gauge transformations of $E^w_\ka$. Define quotient spaces,
$$
\sB^w_\ka=\sA^w_\ka/\sG^w_\ka,\quad
\sB^{w,*}_\ka=\sA^{w,*}_\ka/\sG^w_\ka.
$$
For $A\in A^w_\ka$, let $F_A$ be the curvature of $A$,
let $(F_A^+)_0$ be the  self-dual and
trace-free component of the curvature,
$$
(F_A^+)_0\in \Om^+(\su(E)).
$$
Then define the moduli space of projectively anti-self-dual connections on
$E^w_\ka$ by
\begin{equation}
M^w_\ka=\tM^w_\ka/\sG^w_\ka\quad\text{where}\quad
\tM^w_\ka=\{A\in\sA^w_\ka: (F^+_A)_0=0\}
.
\end{equation}
The following result is well-known.

\begin{thm}
\cite{FU,DK}
For generic choices of metric $g$ on $X$, the moduli space
$M^w_\ka$ is a smooth, orientable manifold of dimension
$$
d(\ka)=8\ka-\frac{3}{2}(\chi(X)+\si(X))
$$
where $\chi(X)$ is the Euler characteristic and $\si(X)$
the signature of $X$.
\end{thm}

\subsection{\Spinu\ structures and Dirac operators}
Let $(W,\rho_W)$ be a \spinc structure on a four-dimensional,
Riemannian manifold $(X,g)$ as defined in
\cite{MorganSWNotes,KrMSWFloer,FL2a}.
  Let $E\to X$ be a rank-two,
complex vector bundle.  A {\em \spinu structure\/} on $X$
is the pair $\ft=(V,\rho)$ where $V=W\otimes E$ and
$\rho:T^*X\to\Hom(V)$ is defined by $\rho=\rho_W\otimes \id_E$.
We define additional bundles $V^\pm=W^\pm\otimes E$ and
$\fg_\ft=\su(E)$.  A more intrinsic definition of these
concepts appears in \cite{FL2a}.  We define characteristic
classes of $\ft$ by
$$
p_1(\ft)=p_1(\fg_\ft), \quad c_1(\ft)=c_1(E)+c_1(W^+),\quad
w_2(\ft)=w_2(\fg_\ft).
$$
A spin connection $A_W$ determines a bijection between
$\sA^w_\ka$ and connections on $V$ respecting the Clifford multiplication
map, $A\mapsto A\otimes A_W$.  This map associates a Dirac operator
$D_A=D_{A\otimes A_W}: L^2_k(V^+)\to L^2_{k-1}(V^-)$
to every $A\in\sA^w_\ka$.
The index of such a Dirac operator is given by
\begin{equation}
\label{eq:DiracIndex}
n_a(\ft)
=
\tquarter\left(p_1(\ft)+c_1(\ft)^2-\si(X) \right).
\end{equation}
\subsection{The degeneracy locus}
For $\ft$ a \spinu structure on $X$ with
characteristic classes $\ka=-\tquarter p_1(\ft)$, $\La=c_1(\ft)$,
and $w$ an integer lift of $w_2(\ft)$, define the
degeneracy locus
\begin{equation}
\label{eq:DefineJ}
J^{\La,w}_\ka
=
\{[A]\in M^w_\ka: \dim\Ker D_A=1\}.
\end{equation}

\begin{thm}
\label{thm:DegenGenMetrics}
\cite{FeehanGenericMetric,TelemanGenericMetric}
Let $\ft$ be a \spinu structure on $X$ with
$n_a(\ft)\le 0$, $\La=c_1(\ft)$, $\ka=-\tquarter p_1(\ft)$,
and $w$ an integer lift of $w_2(\ft)$.
For generic choices of perturbations,
$J^{\La,w}_\ka$ is a smooth submanifold of $M^w_\ka$ of
real codimension $2(1-n_a(\ft))$. The fiber of the normal bundle
of $J^{\La,w}_\ka$ in $M^w_\ka$ at $[A]\in J^{\La,w}_\ka$ is
given by
\begin{equation}
\label{eq:NormalFiber}
\sN_A=\Hom_\CC(\Ker D_A,\Coker D_A).
\end{equation}
Hence,
\begin{equation}
\dim J^{\La,w}_\ka=d_a(\ft)+2n_a(\ft)-2.
\end{equation}
\end{thm}

\begin{defn}
Let $\sN^{\La,w}_\ka\to J^{\La,w}_\ka$ be the
complex rank $1-n_a(\ft)$
normal bundle of $J^{\La,w}_\ka$ in $M^w_\ka$ with
fiber over $[A]\in J^{\La,w}_\ka$ given by $\sN_A$  as defined in
\eqref{eq:NormalFiber}.
\end{defn}

\subsection{The Koschorke-Porteous formula for families of Fredholm operators}
Let $\sF_n(H,H')$ be the space of complex linear, index $n$ Fredholm operators
$L:H\to H'$ where $H$ and $H'$ are Banach spaces.
Let $\sF_{k,n}(H,H')\subset \sF_n(H,H')$ be the subspace of
operators with a $k$-dimensional kernel.  In \cite[\S 3.1]{Koschorke},
Koschorke constructs a cohomology class $\chi_{k,k-n}\in H^{2k(k-n)}(\sF_n(H,H'))$
given, essentially, as the Thom class of the normal bundle of $\sF_{k,n}(H,H')$
in $\sF_n(H,H')$.
Any family, $B$, of complex linear, index $n$ Fredholm operators between separable Hilbert spaces
defines a homotopy class of maps $f:B\to \sF_n(H,H')$ and characteristic classes,
$f^*\chi_{k,k-n}$.  If the family is compact with index bundle $V\in K(B)$,
the characteristic classes $f^*\chi_{k,k-n}$ are related to the Chern classes of $V$
by the equality (see \cite[Thm. 5.2]{Koschorke}),
\begin{equation}
\label{eq:Kosch}
f^*\chi_{1-n,1}=(-1)^{1-n}c_{1-n}(V).
\end{equation}
Further discussions and applications of \eqref{eq:Kosch} appear in \cite[p. 110]{AtiyahJones}.

\section{Chern classes of the index bundle}
\label{sec:DiracChern}
To apply \eqref{eq:Kosch} to prove Theorem \ref{thm:Main}, we
use the following version of the Atiyah-Singer index theorem for families.

\begin{thm}
\label{thm:AS}
\cite{AS4},\cite[Thm. 5.1.16]{DK}
Let $\sE\to X\times B$ be a family of vector bundles with
connections over the
four-manifold $X$ parameterized by $B$.  Let $\fs=(W,\rho)$
be a \spinc structure on $X$.  Let $\DD$ be the index
bundle of the family of Dirac operators,
$$
D_b: \Om^0(W^-\otimes \sE|_{X\times b})
\to \Om^0(W^+\otimes \sE|_{X\times b})
$$
as $b$ varies in $B$.  Then,
$$
\ch(\DD)=-e^{c_1(\fs)/2}\ch(\sE)(1-\frac{1}{24}p_1(X))/[X].
$$
\end{thm}

The first obstruction to applying Theorem \ref{thm:AS}
to the family of Dirac operators $D_A$ for $[A]\in M^w_\ka$
is that the bundle playing the role of $\sE$ for this family,
$$
\EE^w_\ka=\sA^{w,*}_\ka\times_{\sG^w_\ka} E^w_\ka \to \sB^{w,*}_\ka\times X
$$
is not a vector bundle because for $A\in\sA^{w,*}_\ka$,
$\Stab_A=\{\pm 1\}\le \sG^w_\ka$.
As discussed \cite[p.34-35]{RLeungThesis},
the approach to this problem in \cite{RLeungThesis} only
works when the universal $\SO(3)$ bundle defined in \eqref{eq:MuMap}
admits a $\U(2)$ lifting.
To overcome this difficulty in general, we introduce a larger
space with the same rational homotopy type.

\begin{lem}
\label{lem:RationalResolution}
Let $\bS$ be the unit sphere in $L^2(V^+)$.
Then, the map
$$
\pi_S:
\sB^{w,*}_\ka(S)=\sA^{w,*}_\ka\times_{\sG^w_\ka}\bS
\to
\sB^{w,*}_\ka
$$
is a $K(\ZZ_2,1)$-bundle and hence defines an isomorphism
in rational homotopy.
\end{lem}

We define the following subspaces of $\sB^{w,*}_\ka(S)$
analogously,
$$
M^w_\ka(S)=\pi_S^{-1}(M^w_\ka),\quad
J^{\La,w}_\ka(S)=\pi_S^{-1}(J^{\La,w}_\ka).
$$
The following is immediate from Theorem \ref{thm:DegenGenMetrics} and
the surjectivity of $\pi_S$.
\begin{lem}
\label{lem:PullbackNormal}
The space $J^{\La,w}_\ka(S)$ is a smooth submanifold of
$M^w_\ka(S)$ with normal bundle given by
$\pi_S^*\sN^{\La,w}_\ka$.
\end{lem}

For $K^w_\ka\subset M^w_\ka$ a compact subset,
consider the family of Fredholm operators, parameterized
by $(A,\Phi)\in\sA^w_\ka\times\bS$ with $[A]\in K^w_\ka$
\begin{equation}
\label{eq:DiracFamily}
(A,\Phi)\mapsto D_A^*:L^2_k(V^-) \to L^2_{k-1}(V^+).
\end{equation}
This is not a compact family, but it admits a stabilization
in the sense that there is a surjective map from a finite-dimensional, trivial
bundle onto the cokernels of these operators.  The
construction of the stabilization follows immediately from
the independence of
$D_A$ from $\Phi$ and the gauge equivariance of the Dirac
operator (see \cite[Thm. 3.19]{FL2b} for an example of this type
of argument).  This stabilization allows the following definition of
an index bundle for this family of operators.

\begin{defn}
\label{defn:IndexBundle}
For $K^w_\ka\subset M^w_\ka$ a compact subset,
let $\DD^{\La,w}_\ka\in K(\pi_{S}^{-1}(K^w_\ka))$ be the index bundle of the family of
operators defined by \eqref{eq:DiracFamily}.
\end{defn}

Then Lemma \ref{lem:PullbackNormal} and \eqref{eq:Kosch} give us:

\begin{lem}
\label{lem:NormalBundleIndex}
Let $K^w_\ka\subset M^w_\ka$ a compact subset.
Let $\ft$ be a \spinu structure with $n_a(\ft)\le 0$
where $n_a(\ft)$ is the Dirac operator index computed in
\eqref{eq:DiracIndex}.
Then,
$$
\pi_S^*c_{1-n_a(\ft)}(N^{\La,w}_\ka|_{K^w_\ka})
=
c_{1-n_a(\ft)}(\pi_S^*\sN^{\La,w}_\ka)|_{\pi_S^{-1}(K^w_\ka)}
=
(-1)^{1-n_a(\ft)}
c_{1-n_a(\ft)}(\DD^{\La,w}_\ka).
$$
\end{lem}

To apply Theorem \ref{thm:AS} to compute $\ch(\DD^{\La,w}_\ka)$,
we observe that
the index bundle $\DD^{\La,w}_\ka$ is defined as the index bundle
of a family of Dirac operators obtained by twisting
the \spinc structure $(\rho,W^\pm)$ by
\begin{equation}
\label{eq:UniversalU(2)Bundle}
\EE^w_\ka(S)=
\sA^w_\ka\times \bS\times_{\sG^w_\ka} E^w_\ka
\to
\sB^{w,*}_\ka\times X.
\end{equation}
Thus we must compute  $\ch(\EE^w_\ka(S))$.
To this end, we introduce the following cohomology classes.
Recall that for $\beta\in H_\bullet(X;\RR)$, $\mu(\beta)\in H^{4-\bullet}(\sB^{w,*}_\ka;\RR)$
was defined by
\begin{equation}
\label{eq:MuMap}
\mu(\beta)=-\tquarter p_1(\FF^w_\ka)/\beta
\quad\text{where}\quad
\FF^w_\ka=\sA^{w,*}_\ka\times_{\sG^w_\ka}\su(E^w_\ka) \to \sB^{w,*}_\ka\times X.
\end{equation}
If we define
\begin{equation}
\label{eq:UniversalSO(3))Bundle}
\FF^w_\ka(S)
=\pi_S^*\FF^w_\ka
=\sA^w_\ka\times \bS\times_{\sG^w_\ka} \fg_\ft,
\end{equation}
then $\pi_S^*\mu(\beta)=-\tquarter p_1(\FF^w_\ka(S))/\beta$.

\begin{lem}
\label{lem:UnivBundleChar}
Assume $X$ is a smooth, closed four-manifold with $b^1(X)=0$.
Let $E^w_\ka\to X$ be a complex, rank-two vector bundle
with $c_1(E^w_\ka)=w$ and $c_2(E^w_\ka)=\ka+\tquarter w^2$.
Let $\beta_1,\dots,\beta_d$ be a basis for
$H_2(X)/\Tor$ and let $x\in H_0(X)$ be a generator.
Let $\beta^*_i=\PD[\beta_i]$.
Let $Q_{ij}=Q_X(\beta_i,\beta_j)$
and let $P^{ij}$ be the inverse matrix of $Q_{ij}$.
Let $\mu_i=\pi_S^*\mu(\beta_i)$ and $\wp=\pi_S^*\mu(x)$.
Then, for $\EE^w_\ka(S)$ and $\FF^w_\ka(S)$ as defined
in \eqref{eq:UniversalU(2)Bundle} and \eqref{eq:UniversalSO(3))Bundle}
respectively,
\begin{align*}
c_1(\EE^w_\ka(S))
{}&=1\times w,\\
p_1(\FF^w_\ka(S))
{}&=
-4\wp\times 1 -4 \sum_{i,j}P^{ij}\mu_i\times\beta_j^* -4\ka(1\times \PD[x]),
\end{align*}
as elements of rational cohomology.
\end{lem}

\begin{proof}
The first equality follows from observing that
$$
\det(\EE^w_\ka(S))
=
\sA^w_\ka\times \bS\times_{(\sG^w_\ka,\det)} \det(E)
\simeq
\left(\sA^w_\ka\times_{\sG^w_\ka} \bS\right)\times \det(E)
=
\sB^{w,*}_\ka(S)\times\det(E).
$$
The second equality then follows from \cite[5.4.1]{RLeungThesis}
which we now review.
By Lemma \ref{lem:RationalResolution}, $\pi_S$ is an isomorphism in real cohomology
so it suffices to compute $p_1(\FF^w_\ka)$.
By \cite[Prop. 5.1.15]{DK}, $H^\bullet(\sB^{w,*}_\ka;\RR)$ is a polynomial algebra
in $\mu(\beta_i)$ and $\mu(x)$ so we can write
$$
p_1(\FF^w_\ka)
=
a_0\mu(x)\times 1 + \sum_{i,j} a_{i,j} \mu(\beta_i)\times \beta_j^* + b_0 1\times x^*.
$$
To compute the coefficient $a_0$, observe that
$$
p_1(\FF^w_\ka|_{\sB^{w,*}_\ka\times \{x\}})=-4\mu(x).
$$
To compute the coefficient $b_0$, observe that for $[A]\in\sB^{w,*}_\ka$
$$
p_1(\{[A]\}\times X)=p_1(\fg_\ft)=-4\ka\PD[x].
$$
Finally, observe that
\begin{align*}
\mu(\beta_k)
{}&=
-\tquarter p_1(\FF^w_\ka))/\beta_k
\\
{}&=
-\tquarter\sum_{i,j} a_{i,j} \mu(\beta_i)\times \beta_j^*/\beta_k
\\
{}&=
-\tquarter\sum_{i,j} a_{i,j}Q_{j,k}\mu(\beta_i).
\end{align*}
The linear independence of $\mu(\beta_i)$ then implies that
$a_{i,j}=-4 P_{i,j}$ as required.
\end{proof}

Lemma \ref{lem:UnivBundleChar} and
the following will yield the Chern character of $\EE^w_\ka(S)$.

\begin{lem}
\label{lem:ChernChar}
Let $E\to Y$ be a rank two, complex Hermitian vector
bundle.  Let $c_1=c_1(E)$ and $p_1=p_1(\su(E))$.
Then
$$
\ch(E)=
2e^{c_1/2}\sum_{n=1}^\infty \frac{p_1^n}{4^n (2n)!}.
$$
\end{lem}

\begin{proof}
By the splitting principle, we may assume that $E=L_1\oplus L_2$
where $L_i\to Y$ is a complex line bundle.  Let $x=c_1(L_1)$
and $y=c_1(L_2)$.  Then,
\begin{align*}
\ch(E)
{}&=e^x+e^y
\\
{}&=
e^{(x+y)/2}(e^{(x-y)/2} + e^{-(x-y)/2})
\\
{}&=
e^{(x+y)/2}
\left(
    \sum_{n=0}^\infty \frac{(x-y)^n}{2^n n!}
    +
    \sum_{n=0}^\infty \frac{(-1)^n(x-y)^n}{2^n n!}
\right)
\\
{}&=
2
e^{(x+y)/2}
\left(
\sum_{n=0}^\infty \frac{(x-y)^{2n}}{2^{2n}(2n)!}
\right).
\end{align*}
The lemma then follows from the observation that
 $x+y=c_1(E)$ and $(x-y)^2=p_1(\su(E))$.
\end{proof}

\begin{cor}
\label{cor:UniversalChernChar}
Continue the notation of Lemma
\ref{lem:UnivBundleChar}.  Then,
\label{cor:ChernOfUnivBundle}
$$
\ch(\EE^w_\ka(S))
=
2e^{w/2}
\sum_{n=0}^\infty
\frac{(-1)^n}{(2n)!}
\left(
\wp\times 1 +\sum_{i,j} P^{ij}\mu_i\times\beta_j^* +\ka(1\times\PD[x])
\right)^n.
$$
\end{cor}

Our computation of $\ch(\DD^{\La,w}_\ka))$ requires the following algebraic result.

\begin{lem}
\label{lem:UnivP1Products}
Continue the notation of Lemma
\ref{lem:UnivBundleChar}.
Let $h\in H^2(X;\RR)$ satisfy $h=\sum_{k=1}^d h_k \beta_k^*$. Then,
\begin{align*}
(1\times h)\smile (\sum_{i,j} P^{ij}\mu_i\times\beta_j^*)
{}&=\sum_k h_k \mu_k\times\PD[x]
\\
(\sum_{i,j} P^{ij}\mu_i\times\beta_j^*)\smile
(\sum_{k,\ell} P^{k\ell}\mu_k\times\beta_\ell^*)
{}&=
\sum_{i,j}P^{ij}(\mu_i\smile\mu_j)\times\PD[x].
\end{align*}
\end{lem}

\begin{proof}
The first equality follows from
$$
h\smile \beta_j^*=\sum_k h_k\beta_k^*\smile\beta_j^*=\sum_k h_k Q_{kj}\PD[x]
$$
and the equality $\sum_{k} P^{ij}Q_{jk}=\delta^i_k$ where $\delta^i_k$
is the Kronecker delta.  The second equality follows from computing
$$
\beta_j^*\smile\beta_\ell^*=Q_{j\ell}\PD[x],
$$
and the definition of $P^{ij}$ as the inverse of the matrix $Q_{j\ell}$.
\end{proof}

We  now compute $\ch(\DD^{\La,w}_\ka)$.

\begin{prop}
\label{prop:ChernCharOfIndex}
Continue the notation of Lemma
\ref{lem:UnivBundleChar}.
For $c_1(\ft)=\sum_i \la^i \beta_i^*$  define
\begin{equation}
\label{eq:DefineMuTandOm}
\mu(\ft)=\sum_i \la^i\mu_i
\quad\text{and}\quad
\Om=\sum_{i,j} P^{ij}\mu_i\smile\mu_j.
\end{equation}
Then
\begin{equation}
\label{eq:ChernCharOfIndex}
\begin{aligned}
\ch_{2k}(\DD^{\La,w}_\ka)
{}&=
-\frac{(-1)^k}{(2k)!}
\left(
(n_a(\ft)+\frac{2n\ka}{2n+1})\wp^n
-\frac{n}{2(2n+1)}\Om\smile\wp^{n-1}
\right),
\\
\ch_{2k+1}(\DD^{\La,w}_\ka)
{}&=
\frac{(-1)^k }{2(2k+1)!}\ \mu(\ft)\smile\wp^k.
\end{aligned}
\end{equation}
\end{prop}

\begin{proof}
Applying Theorem \ref{thm:AS}, Corollary \ref{cor:UniversalChernChar}
and $c_1(\ft)=c_1(W^+)+c_1(E)$, yields
\begin{equation}
\label{eq:ChernOfIndex1}
\begin{aligned}
{}&
\ch(\DD^{\La,w}_\ka)
\\
{}&
=-\ch(\EE^w_\ka(S))e^{c_1(W^+)/2}\left(1-\frac{1}{24}p_1(X)\right)/[X]
\\
{}&
=
-2e^{c_1(\ft)/2}\left(1-\frac{1}{24}p_1(X)\right)
\sum_{n=0}^\infty
\frac{(-1)^n}{(2n)!}\left(
\wp\times 1 +\sum_{i,j} P^{ij}\mu_i\times\beta_j^* +\ka(1\times\PD[x])
\right)^n/[X].
\end{aligned}
\end{equation}
The second equality of Lemma \ref{lem:UnivP1Products}
implies that
$$
\left(
\sum_{i,j} P^{ij}\mu_i\times\beta_j^* +\ka(1\times\PD[x])
\right)^n
=
\begin{cases}
\sum_{i,j} P^{ij}(\mu_i\smile\mu_j)\times \PD[x] & \text{if $n=2$,}
 \\
0 & \text{if $n>2$.}
\end{cases}
$$
Using the preceding,
we expand the factor in \eqref{eq:ChernOfIndex1}:
\begin{equation}
\label{eq:PowersOfChTerm}
\begin{aligned}
{}&\left(
\wp\times 1 +\sum_{i,j} P^{ij}\mu_i\times\beta_j^* +\ka(1\times\PD[x])
\right)^n
\\
{}&\quad
=
\wp^n\times 1 + n(\wp^{n-1}\times 1)\smile
\left(
\sum_{i,j} P^{ij}\mu_i\times\beta_j^* +\ka(1\times\PD[x])
\right)
\\
{}&\qquad
+
\binom{n}{2}
(\wp^{n-2}\times 1)\smile
\left(\sum_{i,j} P^{ij}\mu_i\times\beta_j^* +\ka(1\times\PD[x])\right)^2
\\
{}&\quad
=
\wp^n\times 1
+
n\sum_{i,j} P^{ij}(\wp^{n-1}\smile\mu_i)\times \beta_j^*
+
n\ka\wp^{n-1}\times\PD[x]
\\
{}&\qquad
+
\binom{n}{2}\sum_{i,j}P^{ij}(\wp^{n-2}\smile \mu_i\smile\mu_j)\times\PD[x].
\end{aligned}
\end{equation}
Then the definition of $\mu(\ft)$ and $\Om$
in \eqref{eq:DefineMuTandOm}, applying
\eqref{eq:PowersOfChTerm}, and the first equality of Lemma \ref{lem:UnivP1Products}
to \eqref{eq:ChernOfIndex1} yield
\begin{align*}
-\ch(\DD^{\La,w}_\ka)
{}&=
2\left(1+\thalf 1\times c_1(\ft) +\frac{1}{8}1\times c_1(\ft)^2 \right)
\left(1-1\times \frac{1}{24}p_1(X)\right)
\\
{}&\quad\smile
\sum_{n=0}^\infty
\frac{(-1)^n}{(2n)!}
\left(
\wp^n\times 1
+
n\sum_{i,j} P^{ij}(\wp^{n-1}\smile\mu_i)\times \beta_j^*
+
n\ka\wp^{n-1}\times\PD[x]
\right.
\\
{}&\qquad\quad
\left. +
\binom{n}{2}\sum_{i,j}P^{ij}(\wp^{n-2}\smile \mu_i\smile\mu_j)\times\PD[x]
\right)/[X]
\\
{}&=
\sum_{n=0}^\infty \frac{(-1)^n}{(2n)!}\wp^n\frac{1}{4}(c_1(\ft)^2-\si(X))
+
\sum_{n=1}^\infty \frac{(-1)^n}{(2n)!} 2n\ka\wp^{n-1}
\\
{}&\quad
+
\sum_{n=1}^\infty \frac{(-1)^n}{(2n)!} n\mu(t)\wp^{n-1}
+
\sum_{n=2}^\infty\frac{(-1)^n}{(2n)!}2\binom{n}{2}\wp^{n-2}\smile \Om
\\
{}&=
\sum_{n=0}^\infty \frac{(-1)^n}{(2n)!}\frac{1}{4}(c_1(\ft)^2-\si(X))\wp^n
-
\sum_{n=0}^\infty
    \frac{(-1)^n}{(2n+2)!} 2(n+1)\ka\wp^{n}
\\
{}&\quad
-
\sum_{n=0}^\infty
    \frac{(-1)^n}{(2n+2)!} (n+1)\mu(t)\smile\wp^{n}
+
\sum_{\ell=1}^\infty
    \frac{(-1)^{\ell+1}}{(2\ell+2)!}2\binom{\ell+1}{2}\wp^{\ell-1}\smile \Om.
\end{align*}
If we observe that $n_a(\ft)+\ka=(c_1(\ft)^2-\si(X))/4$, then
we can write
\begin{align*}
\ch(\DD^{\La,w}_\ka)
{}&=
-\sum_{k=0}^\infty \frac{(-1)^k}{(2k)!}
\left((n_a(\ft)+\ka-\frac{\ka}{2k+1})\wp^k
-\frac{k}{2(2k+1)}\Om\smile\wp^{k-1}
\right)
\\
{}&\quad
-
\sum_{k=0}^\infty \frac{(-1)^k}{2(2k+1)!} \mu(t)\smile\wp^k.
\end{align*}
This completes the proof of Proposition \ref{prop:ChernCharOfIndex}
\end{proof}

\begin{proof}[Proof of Theorem \ref{thm:Main}]
Because $\pi_S^*$ is an isomorphism on rational cohomology,
it suffices to compute the pullback by $\pi_S^*$ of the Chern class.
After Proposition \ref{prop:ChernCharOfIndex}
and Lemma \ref{lem:NormalBundleIndex}, the only remaining step
in the proof of Theorem \ref{thm:Main} is to compute the Chern class
from the Chern character.

If, for a vector bundle $E$, we define
$$
C(E)(t)=\sum_{r\ge 0}c_r(E)t^r,
\quad
Q(E)(t)=\sum_{r\ge 1}q_r(E)t^{r-1},\
\text{where $q_r(E)=r! \ch_r(E)$},
$$
then these classes are related by
\cite[Eq. (2.10')]{MacDonaldSymmetric}
\begin{equation}
\label{eq:GenFunctionRelation}
\frac{d}{dt}C(E)(t)=Q(E)(-t)C(E)(t)\quad\text{or}\quad
C(E)(t)=\exp(\int Q(E)(-t)\ dt).
\end{equation}
The equalities \eqref{eq:GenFunctionRelation} are equivalent to Newton's formula,
\cite[p. 195]{MilnorStasheff} or \cite[(2.11')]{MacDonaldSymmetric},
\begin{equation}
\label{eq:Newton}
c_n
=
-\frac{1}{n}
\sum_{i=1}^n
(-1)^iq_i c_{n-i}.
\end{equation}
From \eqref{eq:ChernCharOfIndex}, we have
\begin{equation}
\label{eq:QClasses}
\begin{aligned}
q_{2k}(\DD^{\La,w}_\ka)
{}&=
(-1)^{k+1}
\left(
(n_a(\ft)+\frac{2k\ka}{2k+1})\wp^k
-\frac{k}{2(2k+1)}\Om\smile\wp^{k-1}
\right)
\\
q_{2k+1}(\DD^{\La,w}_\ka)
{}&=
\frac{(-1)^{k}}{2} \mu(\ft)\smile\wp^k.
\end{aligned}
\end{equation}
Equation \eqref{eq:QClasses} implies that,
abbreviating $\DD=\DD^{\La,w}_\ka$, the power series
$Q(\DD)(t)$ can be written as:
$$
Q(\DD)(t)=\sum_{r=1}^\infty
t^{r-1}
\sum_{i+2j+2k=r} q_{i,2j,2k} \mu(t)^i\smile \Om^j \smile\wp^k,
$$
where, abbreviating $n_a=n_a(\ft)$,
\begin{equation}
\label{eq:NormChernCoeff}
\begin{aligned}
q_{1,0,2k}&= (-1)^{k} 2^{-1},
\quad
q_{0,2,2k-2}  = \frac{(-1)^k k}{2(2k+1)},
\quad
q_{0,0,2k}  =(-1)^{k+1}\left( n_a+\frac{2k\ka}{2k+1}\right).
\end{aligned}
\end{equation}
and $q_{i,2j,2k}=0$ if $i+j>1$.  If we write
$$
C(\DD)(t)=
\sum_{r=0}^\infty \left(\sum_{i+2j+2k=r} f_{i,2j,2k} \mu(t)^i\smile \Om^j \smile\wp^k
\right) t^r,
$$
then the expression for $q_{i,2j,2k}$ in
\eqref{eq:NormChernCoeff} and Newton's formula
\eqref{eq:Newton} imply that the coefficients $f_{i,2j,2k}$ satisfy
the recursion relation:
\begin{equation}
\label{eq:Recursion1}
\begin{aligned}
(i+2j+2k)f_{i,2j,2k}
{}&=
\sum_{u=1}^{k} (-1)^u\left(n_a+\frac{2u\ka}{2u+1}\right) f_{i,2j,2k-2u}
\\
{}&\quad
-
\sum_{u=1}^{k+1} (-1)^u\frac{u}{2(2u+1)}f_{i,2j-2,2k-2u+2}
\\
{}&\quad
+
\sum_{u=0}^{k}\frac{(-1)^u}{2} f_{i-1,2j,2k-2u}.
\end{aligned}
\end{equation}
To find the generating function $C(\DD)(t)$, we compute
$$
P(\DD)(t)=\int Q(\DD)(-t)\ dt
=
\sum_{r=1}^\infty \frac{(-1)^{r-1}}{r}
\left(
\sum_{i+2j+2k=r} q_{i,2j,2k} \mu(t)^i\smile \Om^j \smile\wp^k
\right) t^r.
$$
Then, we write
$P(\DD)(t)=P_1(t)+P_2(t)+P_{3,1}(t)+P_{3,2}(t)$ where
\begin{align*}
P_1(t)&=
\sum_{k=0}^\infty\frac{(-1)^{k}}{2(2k+1)} (\mu\smile \wp^k) t^{2k+1},
\\
P_2(t)&=
\sum_{k=1}^\infty\frac{(-1)^{k+1}}{4(2k+1)}(\Om\smile\wp^{k-1}) t^{2k},
\\
P_{3,1}(t)&=
\sum_{k=1}^\infty\frac{(-1)^kn_a}{2k} \wp^k t^{2k},
\\
P_{3,2}(t)&=
\sum_{k=1}^\infty \frac{(-1)^k\ka}{(2k+1)} \wp^k  t^{2k}.
\end{align*}
If we write
$$
P(\DD)(t)
=
\sum_{r=0}^\infty t^r
\sum_{i+2j+2k=r} m_{i,2j,2k} \mu^i\smile\Om^j\smile\wp^k
$$
and defining
$$
J(x,y,z)= \sum_{i+2j+2k=r} m_{i,2j,2k}x^i y^{2j}z^{2k},
$$
then $J=\thalf xJ_1(z)+\tquarter y^2J_2(z)+J_{3,1}(z)+J_{3,2}(z)$ where
\begin{align*}
J_1(z){}&=\sum_{k=0}^\infty \frac{(-1)^k}{2k+1}z^{2k}=\frac{\tan^{-1}(z)}{z}
\\
J_2(z){}&=
\sum_{k=1}^\infty \frac{(-1)^{k+1}}{2k+1}z^{2k-2}
=
\frac{1}{z^3}\left(-\tan^{-1}(z)+z\right).
\\
J_{3,1}(z)
{}&=
n_a\sum_{k=1}^\infty\frac{(-1)^k}{2k}z^{2k}
=
-\frac{n_a}{2}\ln(1+z^2)
\\
J_{3,2}(z){}&=
\ka\sum_{k=1}^\infty \frac{(-1)^k}{2k+1}z^{2k}
=\ka(\frac{\tan^{-1}(z)}{z}-1).
\end{align*}
This completes the proof of Theorem \ref{thm:Main}.
\end{proof}

\section{Comparing the spin and Donaldson invariants}
\label{sec:Compare}
We now discuss the sources of the error in the equality
\eqref{eq:Incorrect} and the role of Theorem \ref{thm:Main}
in deriving a correct formula for the relation between
the spin and Donaldson invariants.

Define
$$
\AAA(X)=\Sym(H_2(X;\RR)\oplus H_0(X;\RR)).
$$
Then the $\mu$-map of \eqref{eq:MuMap} extends to
$$
\mu:\AAA(X)\to H^\bullet(\sB^{w,*}_\ka;\RR),
\quad
\mu(h^{\delta-2m}x^m)=\mu(h)^{\delta-2m}\smile \mu(x)^m
$$
where $h\in H_2(X;\RR)$ and $x\in H_0(X;\ZZ)$ is a generator.
For $z\in\AAA(X)$,
the definition of a geometric representative, in the sense
of \cite[\S 2.(ii)]{KMStructure}, $\sV(z)$,
for $\mu(z)$ appears in \cite[\S 2.(ii)]{KMStructure}.
If we write $\bar M^w_\ka$ for the Uhlenbeck compactification
of $M^w_\ka$ and $\bar\sV(z)$ for the closure of $\sV(z)$
in this compactification, then the Donaldson invariant
is defined by
$$
D^w_X(h^{\delta-2m}x^m)
=
\begin{cases}
\#(\bar\sV(h^{\delta-2m}x^m)\cap\bar M^w_\ka)
{}&\quad\text{if $\delta=-w^2-3(\chi+\si)\pmod 4$,}
\\
0
{}&\quad\text{otherwise.}
\end{cases}
$$
where $\delta=4\ka-(3/4)(\chi+\si)$.

The error in the equality ,
\begin{equation}
\label{eq:Error}
E^{\La,w}_X(z)=D^w_X(z T^{\La,w}_\ka)-P^{\La,w}_X(z),
\end{equation}
arises from the difference
in the geometric representatives $\sV(T^{\La,w}_\ka)$ and
$J^{\La,w}_\ka$ for the cohomology class $\mu(T^{\La,w}_\ka)$
defined following Theorem \ref{thm:Main}.
The intersections
$$
\sV(T^{\La,w}_\ka)\cap \sV(z)\cap M^w_\ka
\quad\text{and}\quad
J^{\La,w}_\ka\cap \sV(z)\cap M^w_\ka
$$
which define $D^w_X(T^{\La,w}_\ka z)$ and $P^{\La,w}_X(z)$
respectively both have compact support in the interior of a compact subset
$N^w_\ka(z)\subset \sV(z)\cap M^w_\ka$.
While there is a cobordism between the geometric
representatives $\sV(T^{\La,w}_\ka)$ and
$J^{\La,w}_\ka$, this cobordism need not be supported in
$N^w_\ka(z)$. To understand the error \eqref{eq:Error}, one
must therefore study the ends of $\sV(z)\cap M^w_\ka$.

The lower strata of $\bar M^w_\ka$ have the form $M^w_{\ka-\ell}\times\Si$
where $\Si\subset\Sym^\ell(X)$ is a smooth stratum.
Because $\bar\sV(z)$ is, roughly, transverse to the lower strata, the
gluing maps of Taubes, \cite{TauIndef}
(see also \cite{FLKM1} or \cite[\S III.3.4]{FrM}), present the ends of
$\sV(z)\cap M^w_\ka$ as a fiber bundle,
\begin{equation}
\label{eq:FiberBundleEnd}
\begin{CD}
M @>>> \sV(z)\cap M^w_\ka @>>> \sV(z)\cap\left(M^w_{\ka-\ell}\times\Si\right).
\end{CD}
\end{equation}
where the fiber $M$ is a cone, given by the product of moduli spaces of
framed, centered, anti-self-dual connections on $S^4$.

In \cite{RLeungThesis}, Leung computes the error
\eqref{eq:Error}  assuming that $n_a(\ft)= -1$.
He constructs an Ozsvathian ``cap''
(in the sense of \cite{OzsvathBlowUp}) $C^{\La,w}_\ka\subset \sB^{w,*}_\ka$,
for the end of $\sV(z)\cap M^w_\ka$ using the description in
\eqref{eq:FiberBundleEnd}.  If $W^{\La,w}_\ka$
is the intersection of $\sV(z)\cap M^w_\ka$ with the
end \eqref{eq:FiberBundleEnd}, then
$$
\widehat\sV=\left( \sV(z)\cap M^w_\ka - W^{\La,w}_\ka\right)\cup C^{\La,w}_\ka,
$$
will be a compact subspace of $\sB^{w,*}_\ka$ giving the equality
\begin{equation}
\label{eq:CappedEquality}
\#\left(J^{\La,w}_\ka\cap \widehat\sV\right)
=
\#\left(\sV(T^{\La,w}_\ka)\cap\widehat\sV\right).
\end{equation}
Leung then argues that $\sV(T^{\La,w}_\ka)\cap C^{\La,w}_\ka$ is empty
so the right-hand-side of \eqref{eq:CappedEquality} is given by
$D^w_X(z T^{\La,w}_\ka)$ while
the left-hand-side  of \eqref{eq:CappedEquality} is given by
\begin{equation}
\label{eq:LHS}
P^{\La,w}_\ka(z)+\#\left( J^{\La,w}_\ka\cap C^{\La,w}_\ka\right).
\end{equation}
An excision argument (see \cite[Prop. 7.1.32]{DK}
and \cite[\S 4.3 \& \S 6.3]{RLeungThesis})
shows that the intersection number
$$
\#\left( J^{\La,w}_\ka\cap C^{\La,w}_\ka\right)
$$
is independent of the manifold $X$ and can be computed
using examples where both sides of the equality \eqref{eq:Error}
are known.
This argument yields

\begin{thm}
\cite{RLeungThesis}
Continue the hypotheses and notation of Theorem \ref{thm:Main}.
If $n_a(\ft)=-1$, then
\begin{align*}
{}&
D^w_X(h^{\delta-2m-2}x^m T^{\La,w}_\ka) - P^{w,\La}_X(h^{\delta-2m-2}x^m)
\\
{}&\quad
=
\frac{1}{12}
\binom{\delta-2m-2}{2}D^w_X(h^{\delta-2m-2}x^m)Q_X(h)
+\frac{1}{12} mD^w_X(h^{\delta-2m-2}x^{m-1}).
\end{align*}
\end{thm}
Some partial computations for the case $n_a(\ft)=-2$ also appear in
\cite{RLeungThesis}.  However, computing the error \eqref{eq:Error} for
general $n_a(\ft)$ presents some technical challenges.
For $n_a(\ft)<-2$,
the closure $\bar\sV(z)\cap \bar M^w_\ka$ intersects more than one lower stratum
of $\bar M^w_\ka$ and thus more than one open set of the form
\eqref{eq:FiberBundleEnd} is needed to cover the ends of $\sV(z)\cap M^w_\ka$.
Constructing caps in the general case is thus far more
challenging not only because of
the greater topological complexity of the picture \eqref{eq:FiberBundleEnd}
for lower strata but also because the open sets covering the ends overlap.
Some approaches to this type of problem involving multiple open sets have
appeared in \cite{KotschickMorgan,FLMcMaster,FL5,LenessWC}.

Thus, while Theorem \ref{thm:Main} does not give a general relation between
the Donaldson and spin invariants, this result does allow one to localize
the error \eqref{eq:Error} near the lower strata of $\bar M^w_\ka$.
The examples computed by Leung in \cite{RLeungThesis} and similar
computations for the Kotschick-Morgan conjecture, \cite{LenessWC},
and for the $\SO(3)$ monopole cobordism, \cite{FLLevelOne}, suggest
that the error \eqref{eq:Error} has the form described in the following
conjecture.

\begin{conj}
Continue the notation and hypotheses of Theorem \ref{thm:Main}.
The error $E^{\La,w}_X(z)$ of \eqref{eq:Error} is given by a polynomial
in
$$
D^w_X(z O^{\La,w}_\ka),\quad
Q_X,\quad\text{and}\quad
\La,
$$
where $O^{\La,w}_\ka$ is an expression in the characteristic classes
of $\DD^{\La,w}_{\ka-i}$ where $i>0$.
The coefficients of this polynomial depend only on topological
invariants of $X$, $\ka$, and $\La^2$.
\end{conj}.


\providecommand{\bysame}{\leavevmode\hbox to3em{\hrulefill}\thinspace}
\providecommand{\MR}{\relax\ifhmode\unskip\space\fi MR }
\providecommand{\MRhref}[2]{%
  \href{http://www.ams.org/mathscinet-getitem?mr=#1}{#2}
}
\providecommand{\href}[2]{#2}

\bibliographystyle{amsplain}

\begin{thebibliography}{10}

\bibitem{AtiyahJones}
M.~F. Atiyah and J.~D.~S. Jones, \emph{Topological aspects of {Y}ang-{M}ills
  theory}, Comm. Math. Phys. \textbf{61} (1978), 97--118.

\bibitem{AS4}
M.~F. Atiyah and I.~M. Singer, \emph{Index of elliptic operators. {IV}}, Ann.
  Math. \textbf{93} (1971), 119--138.

\bibitem{DonPoly}
S.~K. Donaldson, \emph{Polynomial invariants for smooth four-manifolds},
  Topology \textbf{29} (1990), 257--315.

\bibitem{DK}
S.~K. Donaldson and P.~B. Kronheimer, \emph{The geometry of four-manifolds},
  Oxford Univ. Press, Oxford, 1990.

\bibitem{FeehanGenericMetric}
P.~M.~N. Feehan, \emph{Generic metrics, irreducible rank-one {PU(2)} monopoles,
  and transversality}, Comm. Anal. Geom. \textbf{8} (2000), 905--967,
  math.DG/9809001.

\bibitem{FKLM}
P.~M.~N. Feehan, P.~B. Kronheimer, T.~G. Leness, and T.~S. Mrowka,
  \emph{{PU(2)} monopoles and a conjecture of {M}ari{\~n}o, {M}oore, and
  {P}eradze}, Math. Res. Lett. \textbf{6} (1999), 169--182, math.DG/9812125.

\bibitem{FLKM1}
P.~M.~N. Feehan and T.~G. Leness, \emph{{D}onaldson invariants and
  wall-crossing formulas. {I}: Continuity of gluing and obstruction maps},
  submitted to a print journal, math.DG/9812060 (v3).

\bibitem{FL5}
\bysame, \emph{An {SO}(3)-monopole cobordism formula relating {D}onaldson and
  {S}eiberg-{W}itten invariants}, www.fiu.edu/~lenesst/CobordismBook.pdf.

\bibitem{FL2b}
\bysame, \emph{{$\rm PU(2)$} monopoles and links of top-level
  {S}eiberg-{W}itten moduli spaces}, J. Reine Angew. Math. \textbf{538} (2001),
  57--133, dg-ga/9712005.

\bibitem{FL2a}
\bysame, \emph{{$\rm PU(2)$} monopoles. {II}. {T}op-level {S}eiberg-{W}itten
  moduli spaces and {W}itten's conjecture in low degrees}, J. Reine Angew.
  Math. \textbf{538} (2001), 135--212, math.DG/0007190.

\bibitem{FLLevelOne}
\bysame, \emph{{$\rm SO(3)$} monopoles, level-one {S}eiberg-{W}itten moduli
  spaces, and {W}itten's conjecture in low degrees}, Proceedings of the 1999
  Georgia Topology Conference (Athens, GA), vol. 124, 2002, pp.~221--326.

\bibitem{FLMcMaster}
\bysame, \emph{{$\rm SO(3)$}-monopoles: the overlap problem}, Geometry and
  topology of manifolds, Fields Inst. Commun., vol.~47, Amer. Math. Soc.,
  Providence, RI, 2005, pp.~97--118.

\bibitem{FU}
D.~Freed and K.~K. Uhlenbeck, \emph{Instantons and four-manifolds}, 2nd ed.,
  Springer, New York, 1991.

\bibitem{FrM}
R.~Friedman and J.~W. Morgan, \emph{Smooth four-manifolds and complex
  surfaces}, Springer, Berlin, 1994.

\bibitem{Koschorke}
U.~Koschorke, \emph{Infinite-dimensional {K}-theory and characteristic classes
  of {F}redholm bundle maps}, Global Analysis (F.~E. Browder, ed.), Proc. Symp.
  Pure Math., vol. 15-I, Amer. Math. Soc., Providence, RI, 1970, pp.~95--133.

\bibitem{KotschickMorgan}
D.~Kotschick and J.~W. Morgan, \emph{{SO(3)} invariants for four-manifolds with
  {$b^+=1$}, {II}}, J. Differential Geom. \textbf{39} (1994), 433--456.

\bibitem{KrMSWFloer}
P.~B. Kronheimer and T.~Mrowka, \emph{Monopoles and three-manifolds}, New
  Mathematical Monographs, Cambridge University Press, Cambridge, 2007.

\bibitem{KMStructure}
P.~B. Kronheimer and T.~S. Mrowka, \emph{Embedded surfaces and the structure of
  {D}onaldson's polynomial invariants}, J. Differential Geom. \textbf{43}
  (1995), 573--734.

\bibitem{LenessWC}
T.~G. Leness, \emph{Donaldson wall-crossing formulas via topology}, Forum Math.
  \textbf{11} (1999), 417--457, dg-ga/9603016.

\bibitem{RLeungThesis}
W.-M.~R. Leung, \emph{On $\spinc$ invariants of four-manifolds}, Ph.d. thesis,
  Oxford University, 1995.

\bibitem{MacDonaldSymmetric}
I.~G. Macdonald, \emph{Symmetric functions and {H}all polynomials}, second ed.,
  Oxford Mathematical Monographs, The Clarendon Press Oxford University Press,
  New York, 1995, With contributions by A. Zelevinsky, Oxford Science
  Publications.

\bibitem{MilnorStasheff}
J.~W. Milnor and J.~D. Stasheff, \emph{Characteristic classes}, Princeton Univ.
  Press, Princeton, NJ, 1974.

\bibitem{MorganSWNotes}
J.~W. Morgan, \emph{The {S}eiberg-{W}itten equations and applications to the
  topology of smooth four-manifolds}, Princeton Univ. Press, Princeton, NJ,
  1996.

\bibitem{OzsvathBlowUp}
P.~S. Ozsv{\'a}th, \emph{Some blowup formulas for {SU(2)} {D}onaldson
  polynomials}, J. Differential Geom. \textbf{40} (1994), 411--447.

\bibitem{PTDirac}
V.~Y. Pidstrigatch and A.~N. Tyurin, \emph{Invariants of the smooth structure
  of an algebraic surface arising from the {D}irac operator}, Russian Acad.
  Sci. Izv. Math. \textbf{40} (1993), 267--351.

\bibitem{TauIndef}
C.~H. Taubes, \emph{Self-dual connections on 4-manifolds with indefinite
  intersection matrix}, J. Differential Geom. \textbf{19} (1984), 517--560.

\bibitem{TelemanGenericMetric}
A.~Teleman, \emph{Moduli spaces of {PU(2)}-monopoles}, Asian J. Math.
  \textbf{4} (2000), 391--435, math.DG/9906163.

\end{thebibliography}

\end{document}